\documentclass[a4,12pt,twosided,reqno]{article}
\usepackage[T1]{fontenc}
\usepackage[utf8]{inputenc}
\usepackage{authblk}
\usepackage{lineno,hyperref}
\usepackage{graphicx}
\usepackage{amsmath}
\usepackage{colortbl} 
\usepackage{xcolor} 
\usepackage{algorithm}
\usepackage{algorithmic}
\usepackage[section]{placeins}
\modulolinenumbers[5]
\definecolor{LightCyan}{rgb}{0.80,0.90,1}
\begin{document}
\title{Machine Learning Algorithms in Design Optimization}
\author{Daniele Peri}
\affil{\small CNR-IAC -- National Research Council \\
      Istituto per le Applicazioni del Calcolo "Mauro Picone" \\
      Via dei Taurini 19, 00185 Rome, Italy \\
\texttt{d.peri@iac.cnr.it}}

\maketitle

\begin{abstract} 

Numerical optimization of complex systems benefits from the technological development of computing platforms
in the last twenty years. Unfortunately, this is still not enough, and a large computational time is necessary
for the solution of optimization problems when mathematical models that implement rich (and therefore realistic)
physical models are adopted.

In this paper, we show how the combination of optimization and Artificial Intelligence (AI), in particular
Machine Learning algorithms, can help, strongly reducing the overall computational times, making also
possible the use of complex simulation systems within the optimization cycle. Original approaches are proposed.


\end{abstract}


\section{Introduction}

The very first obstacle to the solution of an optimization problem is represented by the required time. In practice,
the optimization task needs to be included inside the design activities schedule, and a specific time window is
assigned. Consequently, the efficiency of the optimization algorithm is very important, and the number of attempts
(computations of the objective function) before the optimal configuration is identified has to be minimized. But,
although the efficiency of the optimization algorithm is high, the time required for a single evaluation of the
objective function could be so large that some compromises regarding the quality of the physical model become
unavoidable: consequently, only a simplified mathematical model can be applied in practice, and the final solution
is vitiated by this assumption.

For this reason, the use of interpolation/approximation algorithms has been widely adopted, particularly in the
field of Multidisciplinary Design Optimization (MDO), where a cascade of several solvers, one for each discipline
involved, are adopted together, greatly increasing the overall calculation time \cite{Haftka1992}. The first group
of examples is given in \cite{Toropov1992,Barthelemy1992}, more recently other applications can be found
in \cite{Peri2009,Viana2014}. A review can be also found in \cite{Parnianifard2019}. In general, the base idea is
the generation of a {\em meta-model}, that is, {\em a model of the model}, so that an estimate of the objective
function can be obtained by using a simple closed-form expression. A {\em meta-model} can be represented by a
polynomial model \cite{Myers2016,Peri2012b} or something more flexible and
sophisticated \cite{Poggio1990,Matheron1963,Peri2018}.

A great debate about the feasibility of a global interpolation/approximation model for a function of several variables
is still on the table, particularly in the case of a large number of parameters. A sequence of local approximations
has been soon proposed in \cite{Alexandrov1998}, and other techniques in the same framework can be observed
in \cite{Jin2001,Acar2009,Teixeira2019,Peri2012b}: the objective is the reduction of the spatial validity of the
interpolation, localizing the estimate. On the other side, the ambitious goal of building a single
interpolation/approximation model for simulating the state of a system is the typical objective of Artificial Intelligence
(AI), and those techniques, suitable for the management of a great amount of data, can be also adopted for the
description of the digital twin of our physical system.

In this paper, we are describing some techniques able to provide a global interpolation of the state of a system as
a function of the influencing parameters. A further improvement of the {\em meta-model} is then obtained by increasing
the number of samples of the objective function in some critical areas, where the disagreement between the
{\em meta-model} and the true value of the objective function is hypothesized to be high. A regularization technique
is also introduced. The usefulness of these techniques is finally demonstrated thru the application to the solution
of realistic design optimization problems.

\section{Machine Learning for the {\em meta-model} improvement}

As previously noticed, a central point for the optimization of complex systems can be represented by the
determination of a simplified surrogate of a detailed mathematical model of the full system. We are
referring to this as {\em meta-model} because it represents, in practice, {\em a model of the model}.
Here we are recalling some topic elements of the building of the {\em meta-model}.

The first step for the definition of a {\em meta-model} is the generation of a dataset from which we can extract
the information on the optimizing system. This is classically referred to as {\em Design Of Experiments}
(DOE). Since we have typically no information about the function to be fitted, the DOE could be homogeneously
ditributed into the full Design Variable Space (DVS): for the generation of the DOE we can tap into the family of the
so-called {\em Uniformly Distributed Sequences} (UDS). An equally-spaced DOE is also a guaranteee of regularity
of the resulting {\em meta-model}. Examples are the Halton/Hammeresley sequences \cite{Halton1964},
Orthogonal Arrays \cite{Hedayat1999}, Latin Hypercube Sampling \cite{McKay1992}, Sobol sequence\cite{Bratley1988} and
PT-Network \cite{Statnikov1995}.

Furthermore, one may be also interested in evaluating the credibility of the {\em meta-model} in off-design conditions,
that is, the quality of the estimated output when a configuration different from the ones included in the
DOE is checked. Classical theory of Artificial Neural Network (ANN) foresees the division of the DOE in two
subsets, training set and validation set. The training set is used to produce the control parameters of the ANN,
the validation set is used in order to measure the predictive qualities of the ANN.

Since we are considering numerical expensive models (in terms of computational time), the use of a set of data
for validation purposes could appear as a waste of resources. It would be convenient to use all the available
points for the determination of the parameters of the {\em meta-model}, without exceptions. By the way, if a
UDS is adopted, the extraction o a single point is affecting the uniformity of the distribution.

Now we can add some new points for the determination of the performances of the system in some hypothetical
configuration, not included in the training set: this way, we can obtain information about the degree of
precision of the current version of the {\em meta-model}. Once used for their purposes, validation data
should be then added to the training set, re-calibrating the {\em meta-model} with a richer quantity of
information. This is what is commonly called, in AI {\em Machine Learning} (ML), since we are using the
original system (numerical or not) to learn something new to add to our AI system. This could become a
very important phase of the formation of the {\em meta-model} if we could identify a specific area where
it would be useful to add new points to the current {\em meta-model}. This is not easy, since the DOE is
already uniformly distributed on the DVS, so the identification of a new sample cannot be performed based
on purely geometric considerations (i.e., a specific region is not well covered by the training set).

A possible approach to this problem is proposed in \cite{Peri2009}. If we compare different {\em meta-models}
over the full DVS, we observe a different behavior among them, and we cannot determine {\em a priori} which
{\em meta-model} is the best to apply. This situation is typical of a small DOE (undersampling). What we can
do with a moderate computational cost is to compare systematically the outcome of different {\em meta-models}
over the entire DVS, generating a denser UDS for this purpose, and then trace the disagreement between the
prediction of the {\em meta-models}. We can interpret the discrepancy as a measure of the uncertainty in
the prediction so that an additional training point where the disagreement is maximum will surely help in
aligning locally the outputs of the {\em meta-models}. This is far different from adding a new point in an
area with a low density of training points because this approach is explicitly considering the local quality
of the approximation. Numerical experiments reported in \cite{Peri2009} underline improvements of about 10\%
of the quality of the {\em meta-model} for the same number of training points: this is a demonstration that
a uniform distribution of points in the training set is for sure a good start, but customization of the training
set represents a more efficient solution.

\section{Further tuning of the AI model}

A regular distribution of the training points is a prerequisite for the determination of a response surface as
regular. This means that loss of regularity in the distribution of training points can affect the regularity
of the {\em meta-model} response surface, and the previously proposed algorithm is not preserving the regularity
of the distribution: a regularization technique could be helpful in this context.

In the following, we are indicating a possible approach for {\em meta-models} whose construction implies the
solution of a linear system. In particular, we are considering Kriging \cite{Matheron1963,Peri2009} and
Multi-dimensional Spline \cite{Peri2018}.


\subsection{Kriging regularization}

The training of Kriging is performed by assembling and then factorizing the self-correlation matrix of the sample
points \cite{Peri2009}: the spatial correlation between two points of the DVS is determined uniquely based on the
distance between the two points thru the so-called {\em semi-variogram} $\gamma$. In the original formulation,
$\gamma$ is computed experimentally, based on the available dataset. Once the experimental values of the
{\em semi-variogram} are computed, a possible approach is to define $\gamma$ as an exponential function obtained
by fitting the experimental data, whose behavior is typically far from being regular. In theory, there is no reason
why we should not define a different {\em semi-variogram} for each DOE point, so that we are indicating the local
{\em semi-variogram} as

\[
\gamma_i = e^{-(r/a_i)^2}
\]

using different values of $a_i$, one for each point of the DOE. The coefficients $a_i$ are the result of the fitting.
$\gamma_i$ are required for assembling the self-correlation matrix $\Gamma$, whose element $\Gamma_{i,j}$ represents
the correlation between the $i^{th}$ and the $j^{th}$ sample point.  If the coefficients $a_i$ and $a_j$ are different,
$\Gamma$ becomes unsymmetric. In practice, also due to the small amount of DOE points, a single {\em semi-variogram}
is adopted and $\Gamma$ is symmetric. $\Gamma$ is then inverted, and the interpolation of the $N$ DOE values $F(x_i)$
at the generic point $x$ is obtained as

\[
  f(x) = \sum_{i=1}^N w_i F(x_i)
\]

where the weights $w_i$ are obtained by solving

\[
\Gamma W = \Gamma_0
\]

where $\Gamma$ and $\Gamma_0$ are respectievely

\[
\Gamma
=
\begin{bmatrix}
\gamma_1(0)       & \gamma_1(r_{2-1}) & ... & \gamma_1(r_{N-1}) & 1 \\
\gamma_2(r_{1-2}) & \gamma_2(0)       & ... & \gamma_2(r_{N-2}) & 1 \\
 ...              &    ...            & ... & ...               & 1 \\
\gamma_N(r_{1-N}) & \gamma_N(r_{2-N}) & ... & \gamma_N(0)       & 1 \\
 1                &  1                & 1   & 1                 & 0 \\
\end{bmatrix}
\Gamma_0
=
\begin{bmatrix}
\gamma_1(r_{1-0}) \\
\gamma_2(r_{2-0}) \\
 ...              \\
\gamma_N(r_{N-0}) \\
1                 \\
\end{bmatrix}
\]

and $(r_{a-b})$ is $||{\bf x_a}-{\bf x_b}||$. These formulas are used in the Best Linear Unbiased Prediction (BLUP)
of random variables \cite{Robinson}. Under the assumption that the irregularities in the Kriging computation can be
mainly addressed to the condition number of $\Gamma$, we can try to act on the coefficients of the {\em semi-variogram}
$a_i$ in order to maximize the condition number: a compass-search algorithm \cite{Beale1958} is here applied, adjusting
the coefficients $a_i$ maximizing the condition number of $\Gamma$. A maximum variation of $\pm50\%$ is allowed: at
each step, a golden section search \cite{Kiefer1953} is iteratively performed, so that the search limits are easily
enforced.

\begin{figure}[!htb]
\begin{center}
\includegraphics[width=0.95\textwidth]{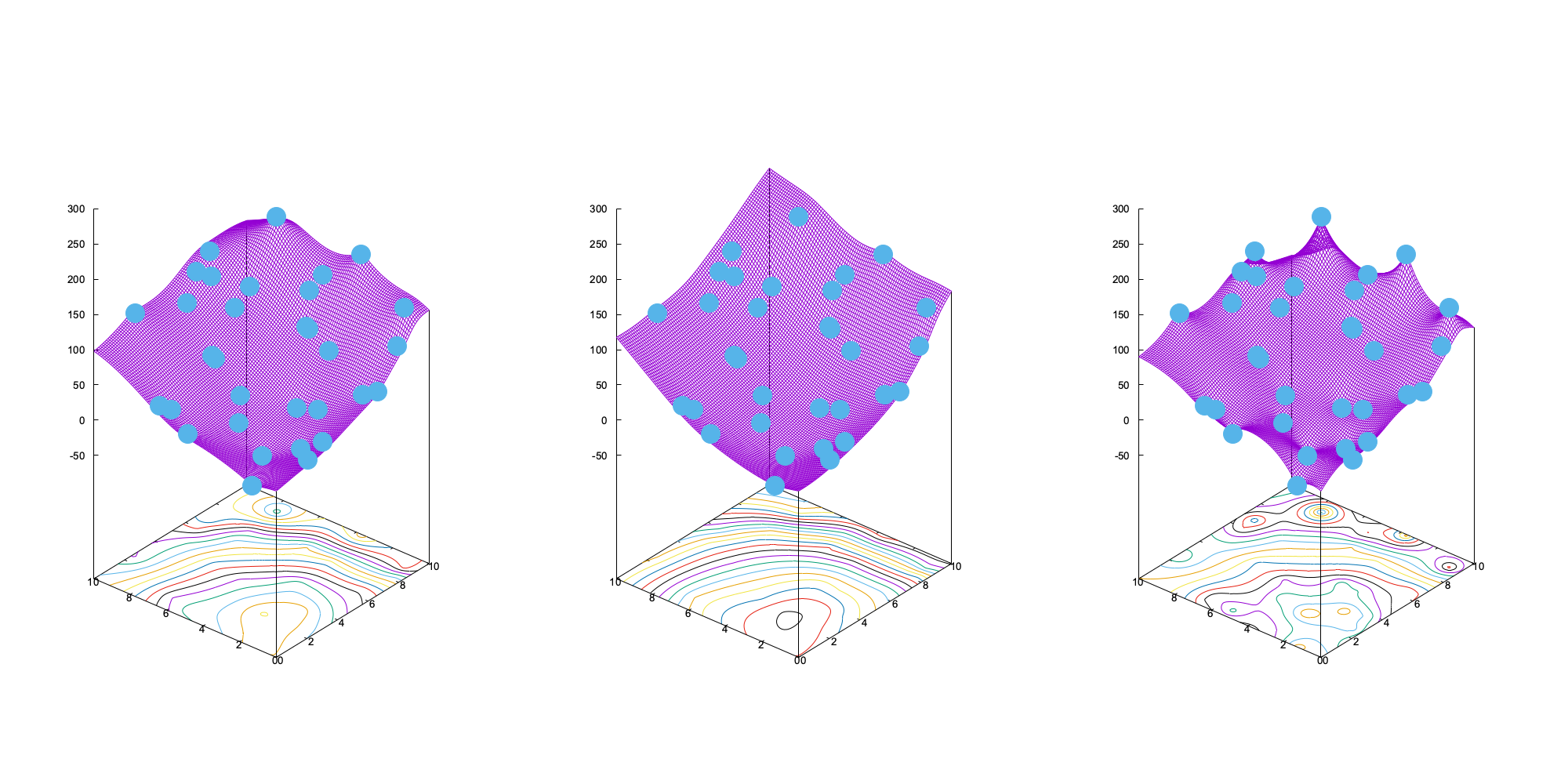}
\caption{Effect of the tuning of Kriging on the base of the condition number of the self-correlation matrix $\Gamma$.
         On the left, the first guess, on the center the effect of a maximization of the condition number (that is, the
         proposed approach), and on right the effect of a minimization of the condition number (that is, the opposite of
         the proposed approach). The test has been performed using the Sasena function in $\Re^2$. The function is sampled
         using 16 random points.
         }\label{evoK}
\end{center}
\end{figure}

In figure \ref{evoK} we can see the effect of the regularization of $\Gamma$ on the overall reconstruction of the objective
function for a 2-dimensional closed form expression (Sasena function
\footnote{$f(x) = 2+0.01*(x_2-x_1^2)^2+(1-x_1)^2+2*(2-x_2)^2+7*sin (0.5*x_1)*sin(0.7*x_1*x_2) \\ x_1 \in [0, 5], x_2 \in [0, 5]$}).
The use of a 2-dimensional function is here justified by the necessity of data visualization. In the first frame on left,
we have the Kriging interpolation, where a single value of $a_i$, obtained by the standard fitting procedure, is applied;
in the central frame, we can see the effect of the maximization of the condition number, with a different value of $a_i$
for each sample point. In the extreme right frame, we have the effect of a minimization of the condition number. The effect
of the maximization of the condition number is evident: the response surface is much more regular, in particular in the
region between the sampled points or in the extreme regions, where there is not a sampled value. On the contrary, a
minimization of the condition number is producing an interpolation formed by several Dirac-type regions, one for each
sample point. It is evident how the criteria appear to be helpful in the regularization of the response surface.

\subsection{Spline regularization}

Also for Multi-Dimensional Spline (MDS), the weights are obtained through the solution of a linear system. Here the
interpolation is obtained as a sum of N compact support functions $R(\rho)$, one for each sample point: also in this
case can, the function can be different from each DOE point:

\begin{equation}\label{eq:pesi}
  f(x) = \sum_{i=1}^N w_i R_i(\rho(x_i,x))
\end{equation}

where $\rho(x_i,x)$ is a measure of the distance between the i$^{th}$ sample point and the computational point $x$,
and $R_i(\rho)$ is a compact support function, decreasing to zero at a certain distance from its center (the sample
point). A simple expression for $R$, linear in the distance between the points, is

\begin{eqnarray*}
 R & = & 1 - b \, ||x_i - x|| for \,\,\, b \, ||x_i-x|| \le 1 \\
 R & = & 0     \,\,\,\,\,\,  otherwise
\end{eqnarray*}

The weights of the kernel functions are determined by solving a linear system enforcing explicitly the equality
between equation \ref{eq:pesi} and the sampled value at every point of the DOE. To compare with the previous case,
while for Kriging the objective function was indirectly included through the {\em semi-variogram}, here it appears
explicitly at the right end side of the linear system.

\begin{figure}[!htb]
	\begin{center}
		\includegraphics[width=0.95\textwidth]{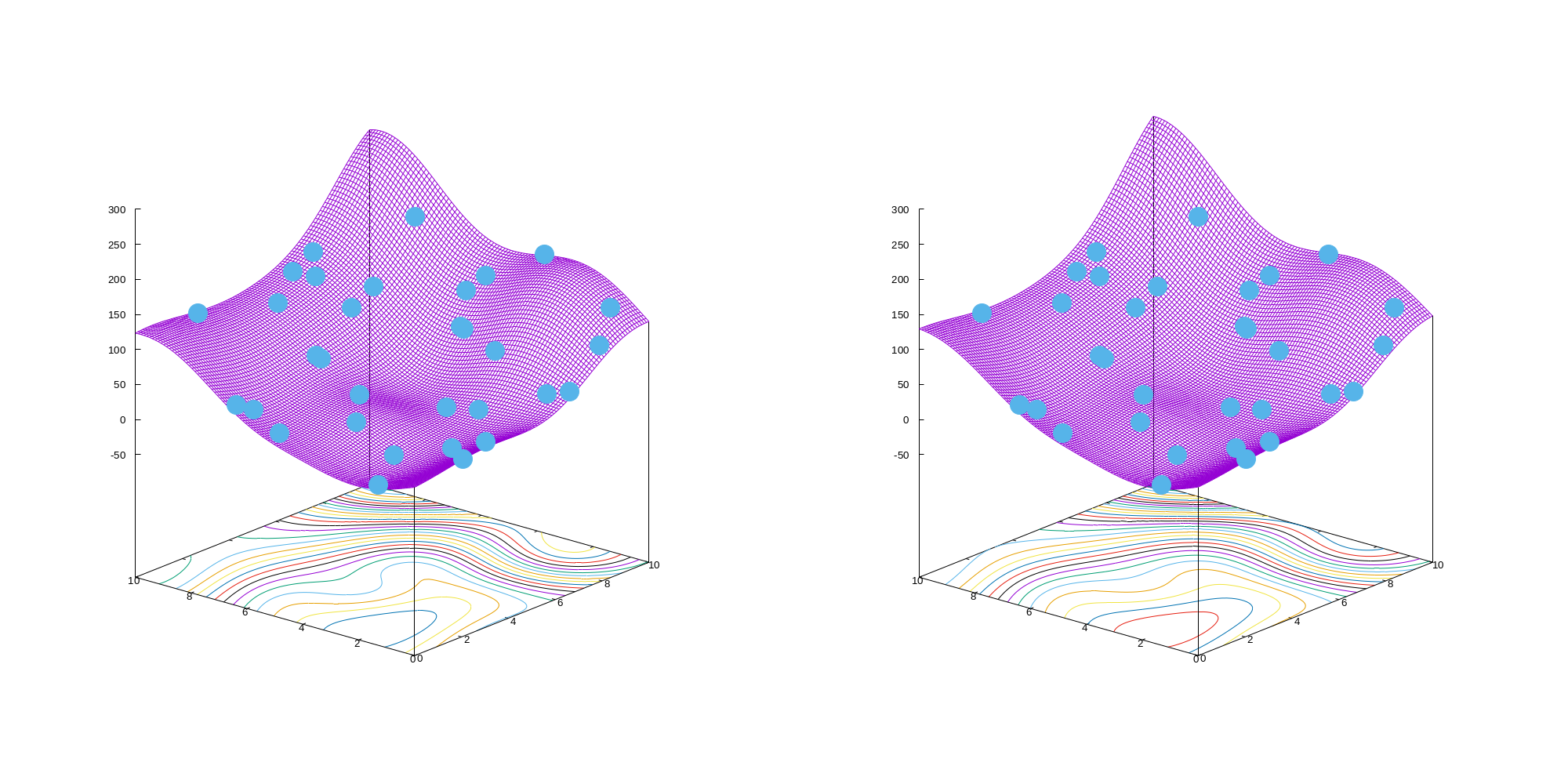}
		\caption{Effect of the tuning of MDS on the base of the condition number of the variance of the
                         coefficients. On the left, is the first guess, and on right is the effect of tuning. The
                         test has been performed using the Sasena function in $\Re^2$. The function is sampled using
                         16 random points.
		}\label{evoN}
	\end{center}
\end{figure}

The parameter $b$ represents a measure of the amplitude of the compact support function, and also this parameter can
be adjusted for the improvement of the condition number of the matrix of the linear system to be solved, using a
procedure equivalent to the one previously described. The effect of the tuning is reported in figure \ref{evoN}.
Compared with Kriging, MDS tends to produce a smoother response surface probably because the weights have a direct
link with the local value of the objective function, and there is a progressive passage between the different influence
areas of the DOE points when we move across the DVS. As observed in the left frame of figure \ref{evoN}, the case
where all the $b_i$s are all equal was already regular. The differences between the first case and the regularized
one are not so large as in the Kriging example, but they can still be observed looking at the contour levels reported
at the bottom of the plot: a more regular behavior is evident in the neighborhood of the minimum of the function. As
a further check, a comparison with the contour lines in figure \ref{evoK} and figure \ref{evoN} indicates that the
effect of the regularization of MDS is also going in the direction of a stricter similitude with Kriging.

\section{Optimization algorithm}

{\em Meta-models} (MM) and ML can be used as base elements for the definition of an optimization algorithm where the
recourse to the mathematical model providing the value of the objective function is very limited. This is a situation
absolutely valuable when the computational cost of a single value of the objective function is very high. Furthermore,
since we are now dealing with a computationally inexpensive surrogate of the objective function, we can also avoid the
recourse to a sophisticated optimization algorithm, proceeding by using a {\em brute force} approach, where the DVS is
sampled extensively by using the MM. The algorithm can be depicted sketchily by the pseudocode reported in
{\bf Algorithm \ref{codice}}.

\begin{algorithm}
\caption{PSI-AI algorithm}\label{codice}
\begin{algorithmic} 

\STATE Perform initial sampling (DOE)

\FOR {$N_E < N_E^{max}$}

\STATE Apply ML 

\STATE Search the DVS using MM 

\STATE Update DOE

\STATE Center the DVS on the current best solution

\STATE Reduce the amplitude of DVS

\ENDFOR
\end{algorithmic}
\end{algorithm}

$N_E$ is the current number of objective function evaluations, and $N_E^{max}$ is the maximum effort we want to put into
the search, measured as the maximum number of evaluations of the objective function we are willing to perform. The name
{\tt PSI-AI} comes from the original formulation of the so-called Parameter Space Investigation (PSI) originally proposed
in \cite{Statnikov1995} but without the ML improvement phase.

The initial sampling of the DVS is performed by using a uniform distribution of points: this is in the logic of the
uniform probability that every point of the DVS could host the optimal value of the objective function unless some
information on the objective function is gained. The DOE is obtained using a P$\tau$-net distribution, extensively
reported in \cite{Statnikov1995}.

\begin{figure}[!htb]
\begin{center}
\includegraphics[width=0.95\textwidth]{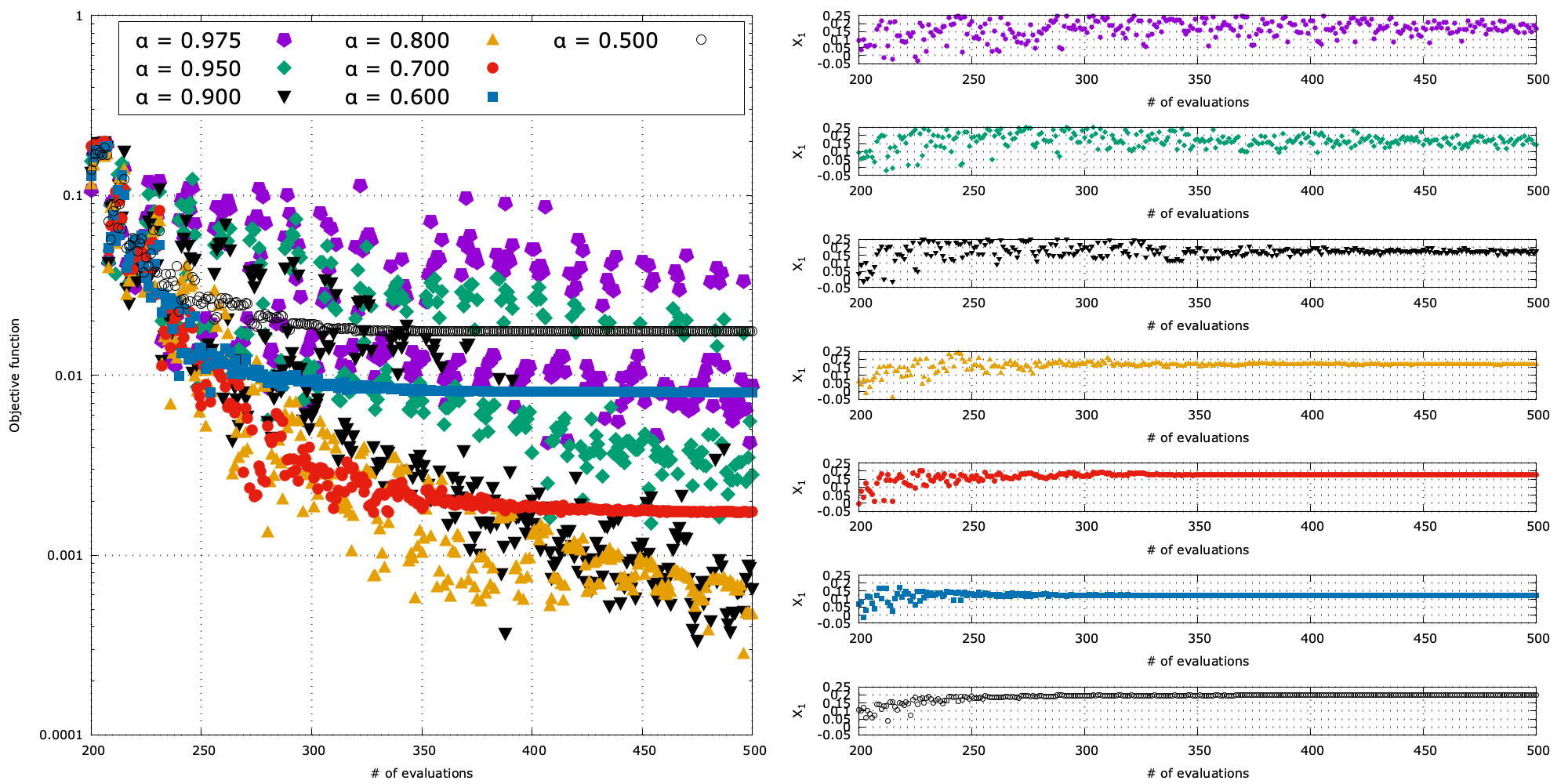}
\caption{Effects of different values of the constriction factor on the convergence of the {\tt PSI-AI} algorithm. On the
         left panel, the percentage difference between the current best value of the objective function and the optimal
         value, and on the right panel the numerical value of the first (of 12) design parameter: different graphs are
         for different values of the constriction factor, with colors in accordance with the left side picture. On the
         horizontal axis, in every graph, the current number of objective function evaluations.
         }\label{xred}
\end{center}
\end{figure}

A relevant parameter is represented by the constriction factor ($< 1$), that is, the amount of reduction of the DVS
amplitude when we are passing from one iteration to another. A small value means a strong reduction: in this case, the
convergence toward the more promising area is fast, but we can lose the location of the global minimum due to premature
focalization. On the contrary, a large value, close to 1, is a guarantee of the completeness of the exploration, but it
can require a very large number of iterations (and then a large number of evaluations of the objective function) to get
the convergence. Here we are performing a comparative study about the role of this parameter. In figure \ref{xred} the
effects of a systematic variation of the constriction parameter $\alpha$ are reported. A quadratic function of 12 variables
is here minimized\footnote{$ f(x) = \sum_{i=1}^{12} (x_i-0.5)^2 \,\,\,\,\, x_i \in [-10:10]$}: the DOE is composed of 192
points ($16 \times$ the number of design variables), and at each iteration, 8 further points are added during the ML phase,
while the best 8 points selected on the base of the evaluation of a regular sampling of the DVS using the MM are also added
at each iteration. 10 iterations are performed. In the picture, to focus on the convergence of the algorithm, the
representation of DOE evaluation is not included.  We can observe from figure \ref{xred} that a value of $\alpha=0.5$ causes
a premature convergence, and the reduction of the objective function is quite large but also far from the best value. A
value of $\alpha$ in between 0.8 and 0.9 appears to be the most appropriate choice, with a slight preference for 0.9,
since it looks like a further improvement of the objective function would be possible if a larger number of iterations
were allowed. Larger values appear to slow down too much the convergence of the algorithm, although they are a guarantee
(asymptotically) of a more accurate exploration of the DVS.

\section{Application to realistic problems}

\begin{figure}[!htb]
\begin{center}
\includegraphics[width=0.95\textwidth]{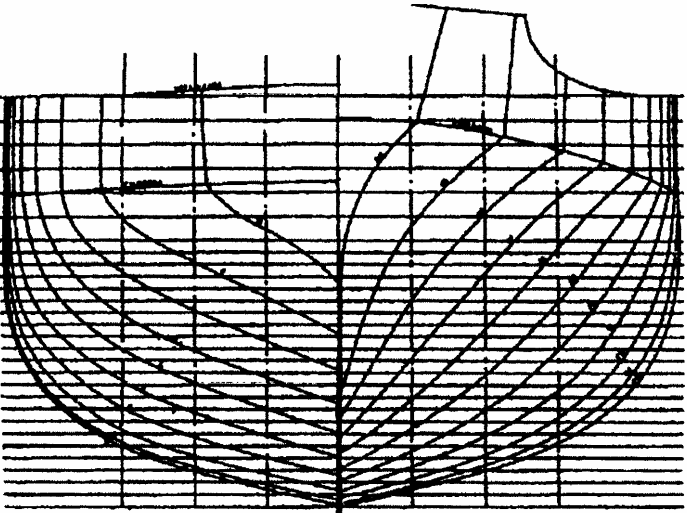} \\
\begin{tabular}{c|c|c|c}
L$_{OA}$  & 22.98 [m] & Maximum beam      & 4.22 [m]       \\ \hline
L$_{PP}$  & 21.85 [m] & Beam at waterline & 4.20 [m]       \\ \hline
Draught   &  1.40 [m] & Displacement      & 53.15 [t]      \\ \hline
          &           & Wetted surface    & 83.35 [m$^2$]  \\
\end{tabular}
\caption{Section lines of the {\em Vaproetto} as from \cite{Zotti2002}.
         }\label{Vaporetto}
\end{center}
\end{figure}

To check the qualities of the algorithm on a realistic application, whose resulting objective function is possibly not
showing a single clear global minimum as in the adopted algebraic test case, the problem of the optimization of a monohull
ship has been considered. The parent hull form (PHF) for the ship design optimization has been taken from \cite{Zotti2002}:
it represents the bare hull geometry of a {\em Vaporetto}, the water-bus performing public transport in Venice. Lines and
main dimensions are reported in figure \ref{Vaporetto}.

\subsection{Pitch motion reduction}

As a first case study, the {\tt PSI-AI} algorithm has been compared with an efficient multi-agent optimization algorithm,
to understand and evaluate the performance we can expect from {\tt PSI-AI}. The selected optimization algorithm is the
hybrid version of the Imperialist Competitive Algorithm ({\tt ICA}), originally proposed in \cite{Atashpaz2007} and then
here adopted in the improved version \cite{Peri2019} as {\tt hICA}. The {\tt hICA} includes also a local search algorithm,
used in conjunction and in cooperation with the original multi-agent algorithm.  {\tt hICA} has been proved to be more
efficient than other multi-agent algorithms, like the NSGA-II implementation \cite{Deb2002} of a Genetic Algorithm
({\tt GA})\cite{Peri2019}.

To limit the computational effort, the energy associated with the pitch response of the ship in rough seas has been
considered as objective function. This choice allows the application of a very fast (but reliable) simulation tool,
the open-source {\tt PDStrip} seakeeping code \cite{Bertram2006}. The objective function is represented by the area
below the RAO pitch curve, at the speed of 10 knots. Its evaluation has a computational cost lower than 9 seconds on
a 3.30GHz Intel\textregistered Core\texttrademark i7-5775R: as a consequence, we can easily perform extensive tests
in a reasonable time.

\begin{figure}[!htb]
\begin{center}
\includegraphics[width=0.95\textwidth]{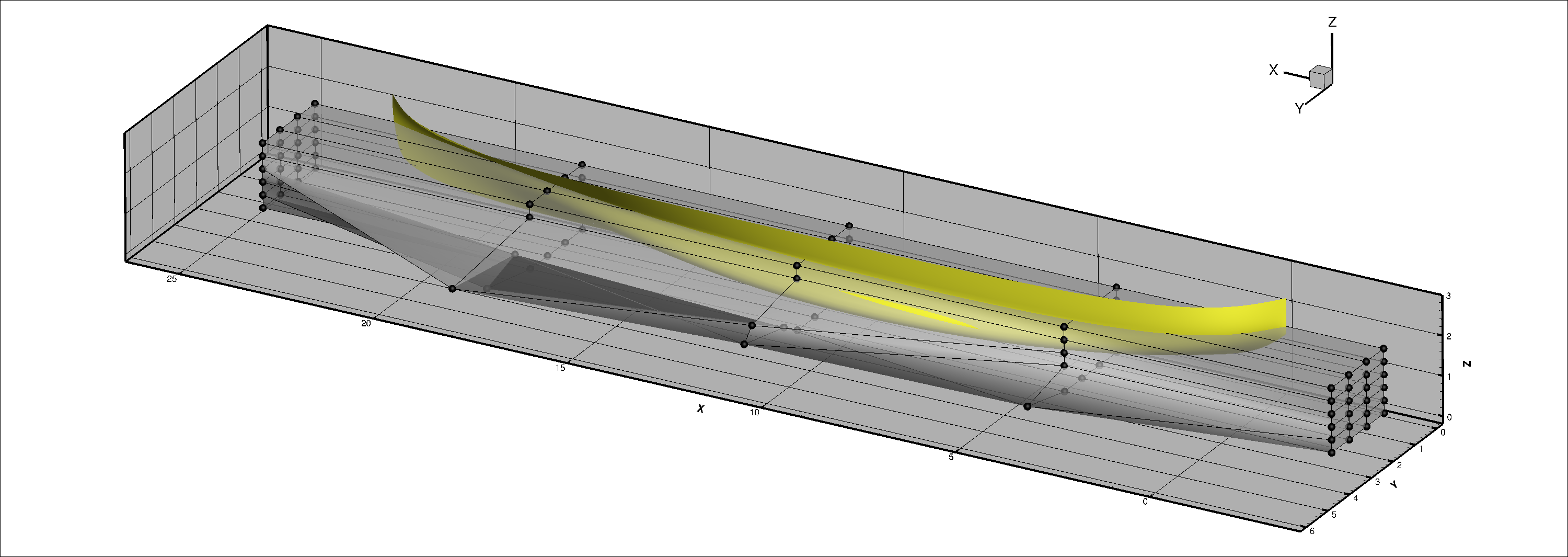}
\caption{Example of shape deformation of the {\em Vaporetto} geometry. The hull is partly embedded into a single FFD patch,
         leaving unchanged the top part of the ship. Here the deformation producing the optimal shape is reported. Black
         dots are the control points of the FFD, while the hull shape is represented in yellow. Bow points toward the
         positive direction so that we are here observing the hull from the stern.
         }\label{FFD}
\end{center}
\end{figure}

The parameterization of the hull has been produced by using the Free Form Deformation (FFD) approach, proposed and
described in \cite{Sederberg1986}. A patch with $5 \times 4 \times 6$ subdivisions (respectively along the {\bf X},
{\bf Y} and {\bf Z} Cartesian axes) has been adopted, but only 12 control points are active: the first and the last
planes along the {\bf X} direction are kept fixed, as well as three first slices in the {\bf Y} direction together
the bottom plane and the first two top planes in the {\bf Z} direction. Only the $4 \times 3$ control points on the
last lateral slice parallel to the {\bf XZ} plane can freely move along the lateral direction, their movement is limited
to 25 centimeters (ship scale). The patch is including the PHF geometry up to the plane $z$ = 1.5 meters, in order not
to change the bridge geometry: the fixed control points of the FFD enforce the continuity between the modified and
unmodified parts of the hull. An example of possible deformation obtained with this configuration is reported in
figure \ref{FFD}.

Among the different possible choices, Kriging has been adopted as a substitute for the mathematical model, while the MDS
has been adopted for the ML phase together with Kriging.

\begin{figure}[!htb]
\begin{center}
\includegraphics[width=0.95\textwidth]{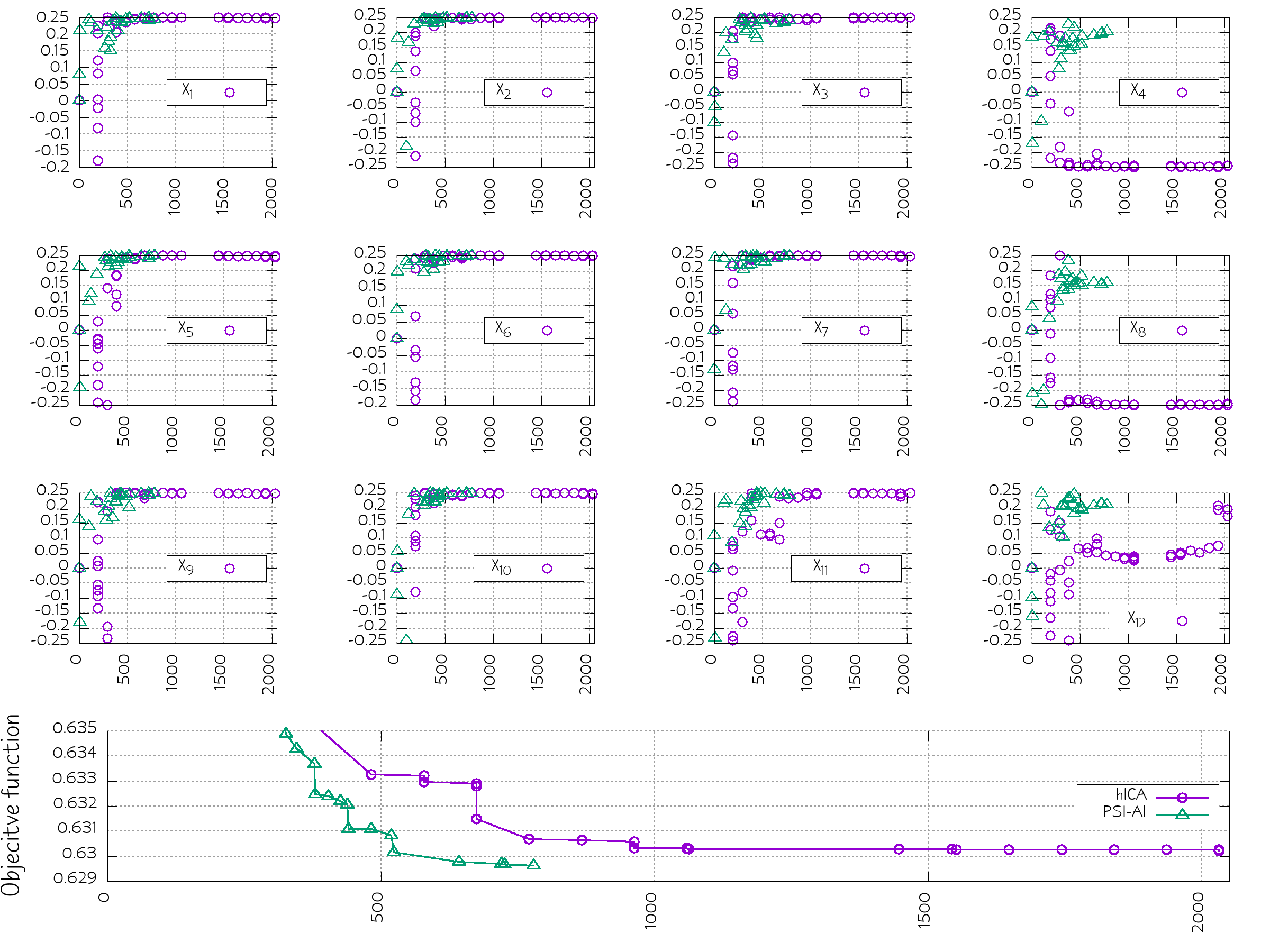} \\
\caption{Comparison between the {\tt PSI-AI} approach and the {\tt hICA} algorithm. On the top, is the time history of
         the 12 design variables. On the bottom, the progressive reduction of the objective function value thru the
         iterations. A dot is plotted only if an improvement of the objective function occurs.
         }\label{bestRAO}
\end{center}
\end{figure}

The results, in terms of the objective function and design parameters, are reported in figure \ref{bestRAO}.  Here is
clear how the number of objective function evaluations required to reach convergence for {\tt PSI-AI} is about one-third
with respect to {\tt hICA}. On the other hand, we have some differences in the optimal values of the design parameters,
whose numerical values are reported in table \ref{tab:para}.

Firstly we analyze the evolution of the design parameter values thru the optimization procedure.  We can clearly observe
how the optimal values of the design parameters in many cases are similar, but not in all the cases. Since most of the
parameters find their optimum value on the border of the admissible portion of the DVS, we can argue that there is not
a stationary point in the constrained DVS, so the selected optimum is the one producing the best possible value of the
objective function, but the local derivatives of the objective function are not null. This situation is clear for 9 out
of 12 parameters.  For the remaining 3 parameters (\#4, 8, and 12), the difference is mainly addressable to their small
sensitivity of the objective function: in fact, the final value of the objective function for the two optimal points is
absolutely comparable. The two algorithms have different attitudes in the exploration of the DVS. In {\tt PSI-AI} the
search has not a preferential area, and the focusing is progressive, and the initial search is mainly grounded on the
center of the DVS. On the contrary, the {\tt hICA} generally expands towards the extreme regions of the admissible DVS,
eventually resizing the group of the agents in a second moment. As a consequence, if the sensitivity of a parameter is
small, {\tt PSI-AI} will preferentially return a value not close to the borders of the admissible DVS, while {\tt hICA}
is implicitly encouraging the values at the extremes.

\begin{table}[!htb]
\caption{Optimal values of parameters obtained using the two different optimization algorithms, {\tt PSI-AI}
        and {\tt hICA}. The two sets of parameters are represented side-by-side, ordered by the
        number of the design parameter.
         }\label{tab:para}\centering\small
\begin{tabular}{l|l|l||l|l|l||l|l|l}
\# & {\tt hICA}  & {\tt PSI-AI}   & \# & {\tt hICA}  & {\tt PSI-AI} & \# & {\tt hICA}  & {\tt PSI-AI} \\ \hline
 1 &  0.2490     & 0.2432         &  2 &  0.2484     & 0.2498       &  3 &  0.2495     & 0.2499       \\ \hline
\cellcolor{orange}{4}             & 
\cellcolor{orange}{-0.2454}       & 
\cellcolor{orange}{0.2045}        & 
 5 &  0.2441     & 0.2498         &  6 &  0.2499     & 0.2500       \\ \hline
 7 &  0.2486     & 0.2496         &
\cellcolor{orange}{8}             & 
\cellcolor{orange}{-0.2494}       & 
\cellcolor{orange}{0.1592}        & 
 9 &  0.2490     & 0.2484	         \\ \hline
10 &  0.2455     & 0.2485         & 11 &  0.2493     & 0.2412       &
\cellcolor{LightCyan}{12}            & 
\cellcolor{LightCyan}{ 0.1957}       & 
\cellcolor{LightCyan}{ 0.2105}       \\ 
\end{tabular}
\end{table}

\begin{figure}[!htb]
\begin{center}
\includegraphics[width=0.95\textwidth]{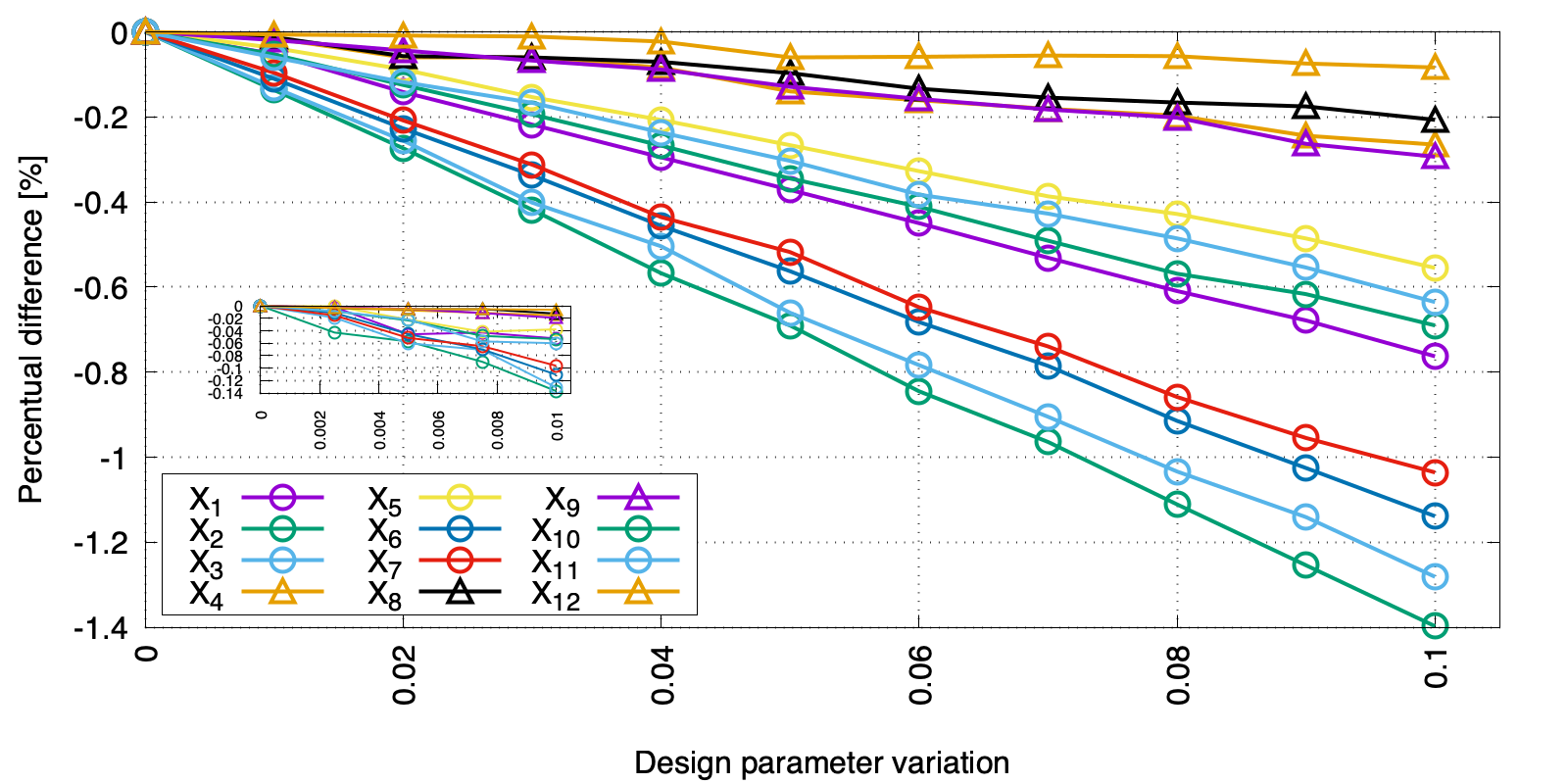}
\caption{Sensitivity analysis of the pitch peak as a function of the design parameters. In the sub-picture, the detail
         of the effect of smaller variations of the design parameters: here the horizontal scale is about one-tenth of
         the larger picture.
         }\label{sensi}
\end{center}
\end{figure}

This hypothesis is substantially confirmed by the results of the sensitivity analysis reported in figure \ref{sensi}.
In this picture, a single parameter is changed at a time, revealing the influence of each parameter on the objective
function. Here we can observe how there is a small group of design variables with a reduced sensitivity with respect
to the objective function: all the previously listed parameters belong to this group.  As a consequence, the small
influence of these parameters is the reason why they are much more prone to follow the tendency subtly suggested by
the algorithm.

\begin{figure}[!htb]
\begin{center}
\includegraphics[width=0.95\textwidth]{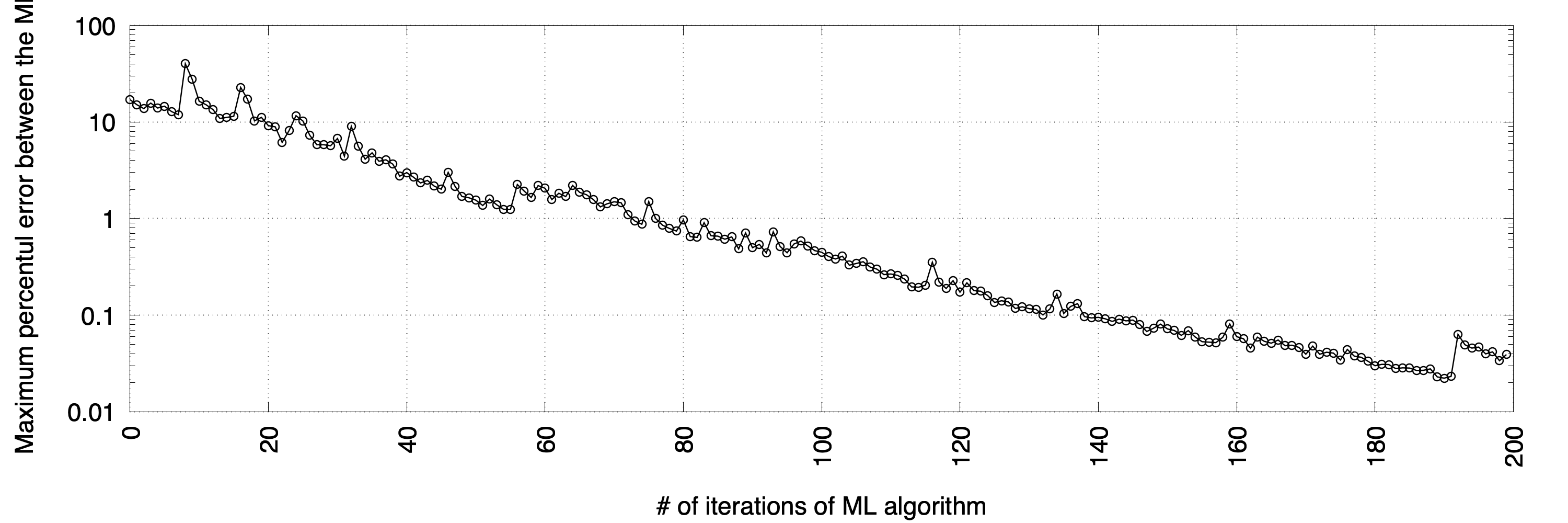}
\caption{Effect of the Machine Learning algorithm on the relative precision of the two different {\em meta-models}.
         On the horizontal axis, the number of calls to the objective function required by the ML algorithm during
         the optimization problem, on the vertical axis the maximum value across the full DVS of the percentage
         difference in the prediction of the two different {\em meta-models}. A reduction of the discrepancies is
         a guarantee of the reliability of the interpolations.
         }\label{MachineLearning}
\end{center}
\end{figure}

In figure \ref{MachineLearning} the progressive reduction of the discrepancies between the two adopted MMs during the
course of the optimization process is reported. The full number of calls to the objective function is about 800: this
means that only a quarter of them are devoted to the ML algorithm. After 80 evaluations, the difference is of the
order of 1\%: this can be considered a relevant achievement, demonstrating the efficiency of the ML algorithm.

\subsection{Total resistance}

A second test case has been produced to further check the capabilities of the {\tt PSI-AI} algorithm. The objective
function is now represented by the effective power in calm water of the {\em Vaporetto} at the speed of 4.5 m/s, which
is a little lower than the maximum speed fixed by the rules (20 km/h). Due to the motivations of the present work,
the shape of the sea bottom and the side walls of the channels have been not considered, although they are peculiar
elements for this kind of sea vehicle traveling in the Venice area.  Hydrodynamic computations are performed using a
single-layer potential flow solver, of the class of the Boundary Element Methods (BEM) \cite{Gadd1976}. These kinds
of solvers are correctly modeling the wave pattern generated by the hull, but the viscous effects are not included
in the formulation, so they are introduced {\em a posteriori} in a simplified way. The resulting estimate of the
effective power has been proved to be reliable, at least from the engineering point of view. The computational effort
of a fully-3D BEM is larger than a strip-theory method, also because computations need to be repeated iteratively to
obtain the actual values of the sinkage and trim of the ship. The CPU time for a single value of the objective function
is now about 40 seconds on the same computational platform. For this reason, the comparison with the {\tt hICA} algorithm
has been not repeated.

The same parameterization as from the previous test case has been adopted, including also the constraints on the design
variables.

At the end of the optimization problem, the effective power required by the {\em Vaporetto} is passing from 22.44 to 15.10 kW,
with a reduction of 32.72\%. To have a second check of the improvements, the RANSE solver {\tt interFoam}, from the suite
{\tt OpenFOAM}\cite{OpenFOAM} has been also applied: results are reported in table \ref{tab:ave}, where we can see the
substantial equivalence of the estimated percentual reductions (34.1 vs. 32.7\%). In this case, the Reynolds-averaged
Navier-Stokes equations, including explicitly the viscous terms neglected in the Laplace equation of the BEM method,
have been solved. Also {\tt interFoam} is taking into account the real sinkage and trim of the hull.

\begin{table}[!htb]
\caption{Optimal values of the wave resistance, frictional resistance, and total resistance for the original and
         optimized hulls, computed by using the RANSE solver {\tt interFoam} from the suite
         {\tt OpenFOAM}\textcopyright \cite{OpenFOAM}. The simulations are performed with the hull able to take the
         dynamic sinkage and trim. Since the effective power differs by the total resistance by a constant, the
         percentage differences of total resistance are the same as the ones of the effective power.
         }\label{tab:ave}\centering\small
\begin{tabular}{l|c|c|c}
Ship      & R$_W$ [N] & R$_F$ [N] & R$_T$ [N] \\ \hline
Original  & 2411.0    & 430.9     & 2841.9 \\ \hline
Optimal   & 1638.9    & 234.7     & 1873.6 \\ \hline\hline
\cellcolor{orange}{$\Delta$\%} & 
\cellcolor{orange}{-32.0} &
\cellcolor{orange}{-45.5} &
\cellcolor{orange}{-34.1} \\
\end{tabular}
\end{table}

\begin{figure}[!htb]
\begin{center}
\includegraphics[width=0.95\textwidth]{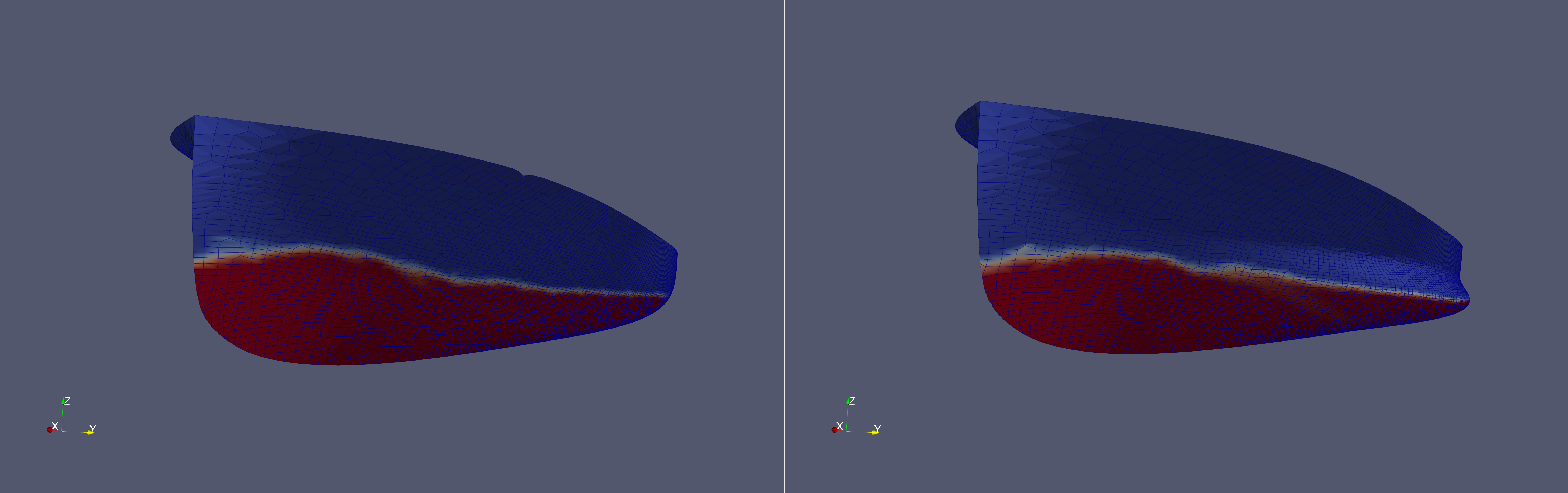}
\caption{Perspective view of the PHF (on left) and the optimal (on right) hull shapes. The red color is
         indicating the wetted part of the hull (as computed by using the {\tt interFoam} RANSE solver.
         }\label{geometry}
\end{center}
\end{figure}

Figure \ref{geometry} is reporting the changes in the shape of the PHF. The beam is increased, and a part of the volume
is shifted fore: since the displacement is fixed, the draught of the optimal hull is smaller in comparison with the PHF
one. The topside of the hull is not changing: as a consequence, the beam increase is causing a kind of lateral bulb in
the central part of the ship. If required, the hull lines of the out-of-the-water part of the hull, not directly
influencing the performance, can be reassessed to have a more regular shape. The red color on the hull geometry is
evidencing the wetted part of the hull at the optimizing speed.

\begin{figure}[!htb]
\begin{center}
\includegraphics[width=0.95\textwidth]{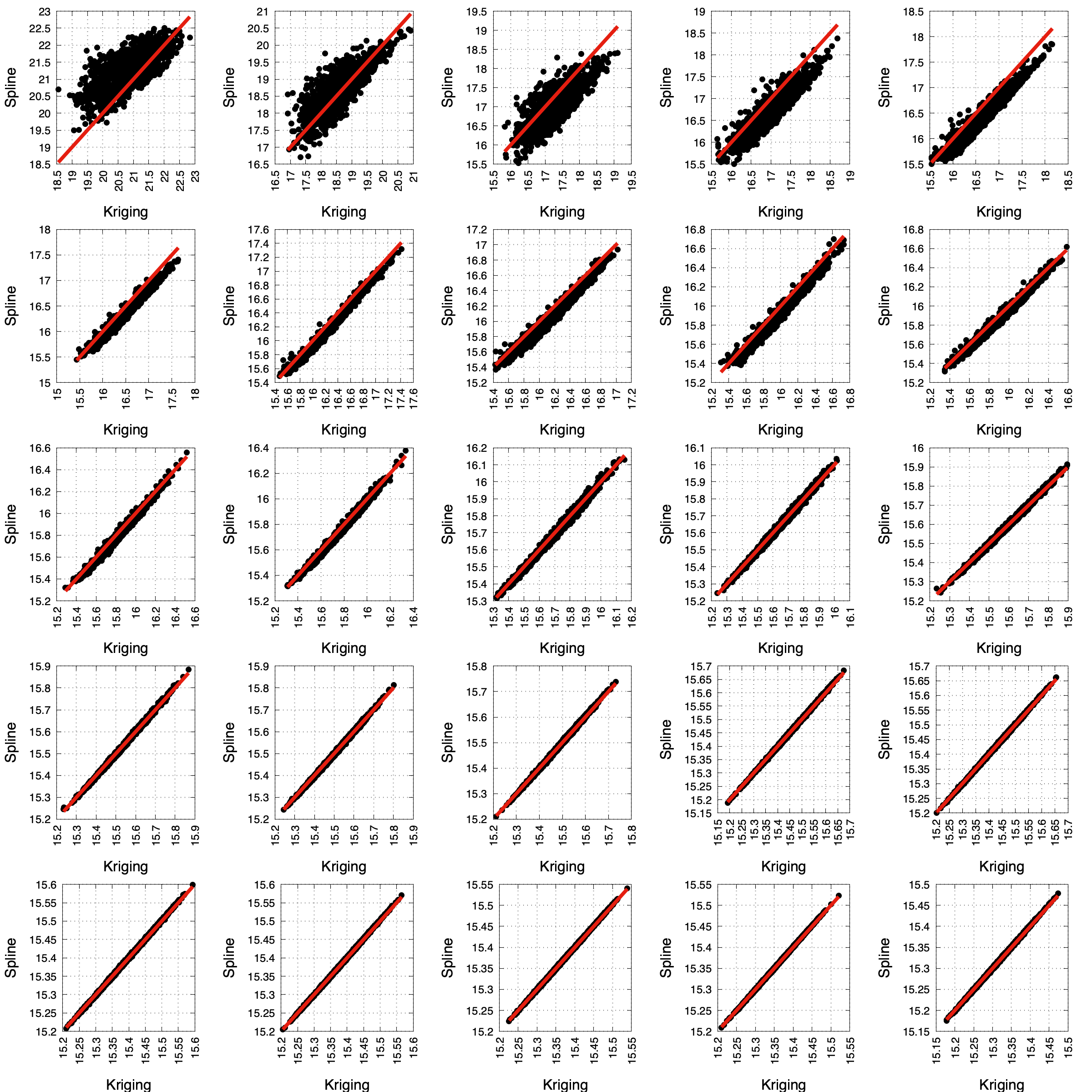}
\caption{Effect of the application of the ML algorithm on the overall precision
         of the adopted {\em meta-models}. In each framework, the values
         of the objective function as predicted by the two different {\em meta-models} 
         at a specific iteration are reported on the two axes. Iteration number (frame number)
         is running from left to right, top to bottom.
         }\label{correlation}
\end{center}
\end{figure}

In figure \ref{correlation} we have reported a different way to stress the increasing similitude between the prediction
of the two MMs as the ML iterations are going on. The two MMs have been computed on a large number of points
laying into the DVS: on each frame, we have on the horizontal axis the estimate provided by Kriging, while the estimate by
MDS is reported on the vertical axis. If the values estimated by two MMs were the same, all the points would be aligned
along the $x=y$ line in the plot. If the points are not well aligned, there is still a difference between the prediction
of the two MMs. We can see how, proceeding from left to right, top to bottom, the points are thickening on the line of
full correlation: this is a sign of a progressive increase in the coherence between the two MMs. Since all the points are
aligned, the similitude of the two MMs is certified on the full DVS.

\begin{figure}[!htb]
\begin{center}
\includegraphics[width=0.95\textwidth]{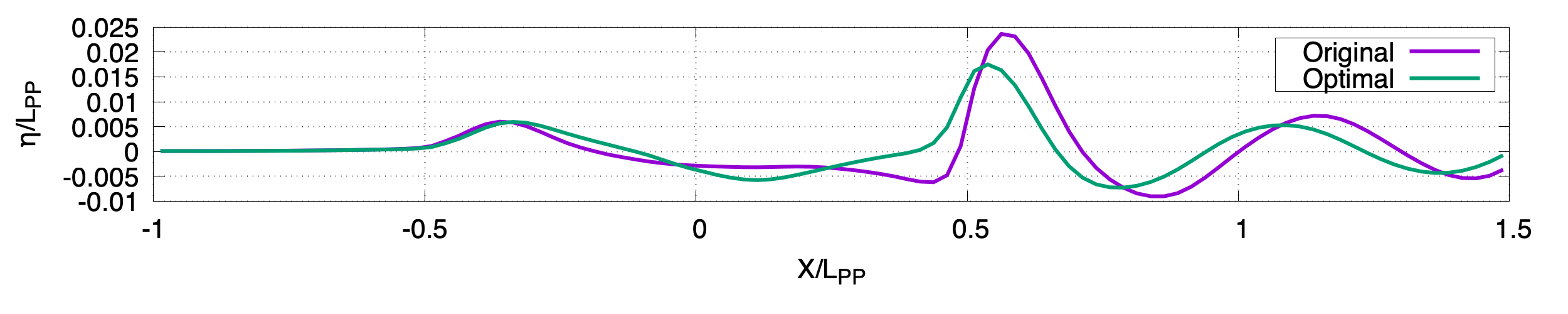}
\caption{Comparison between wave pattern generated by the PHF and optimal hulls,
         computed with the BEM solver \cite{Gadd1976}, the one applied for the evaluation
         of the objective function. The ship hull is placed in between $x=-0.5$ and $x=0.5$,
         with bow at $x=-0.5$.
         }\label{etaBEM}
\end{center}
\end{figure}

\begin{figure}[!htb]
\begin{center}
\includegraphics[width=0.95\textwidth]{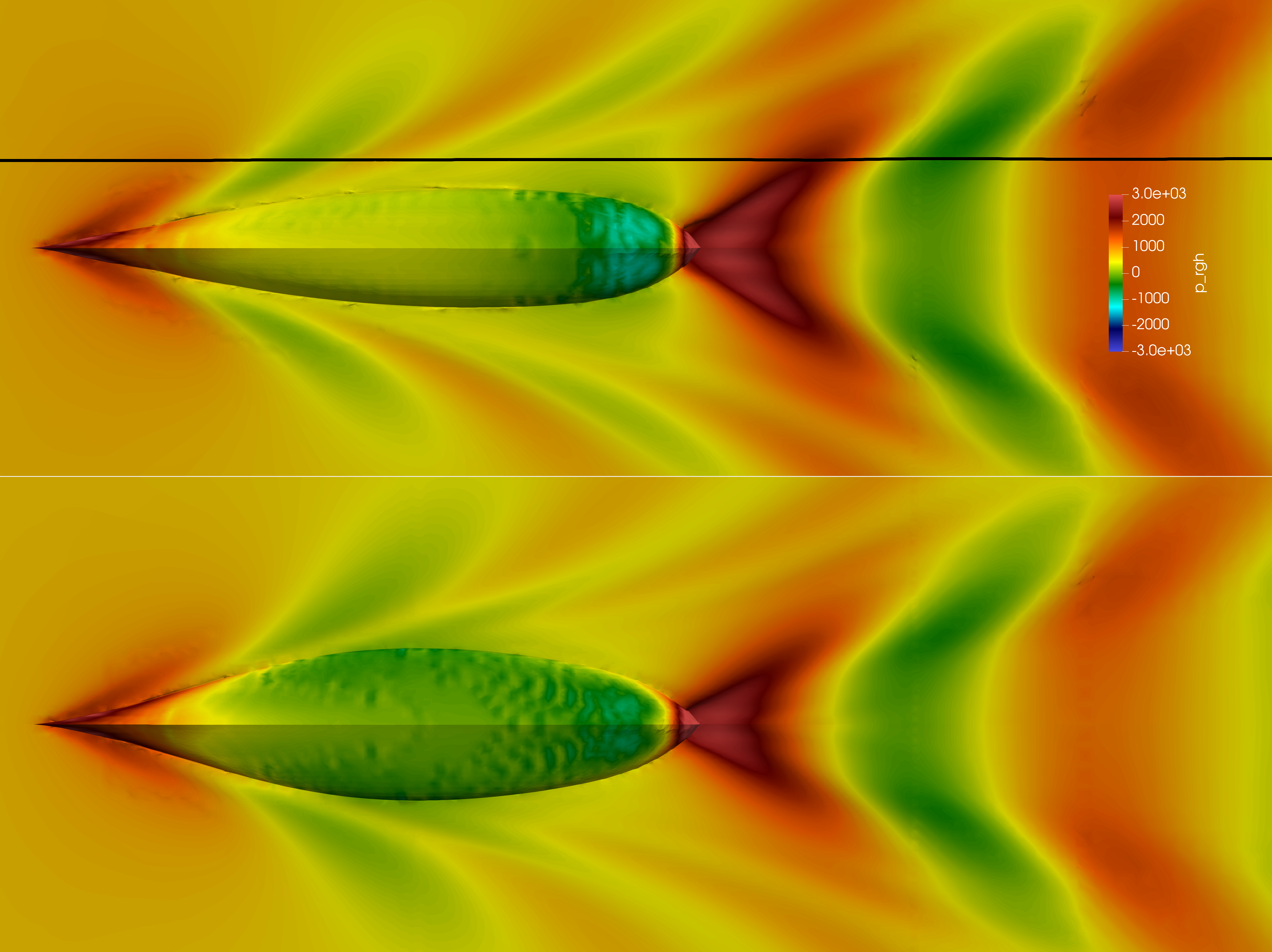} \\
\includegraphics[width=0.95\textwidth]{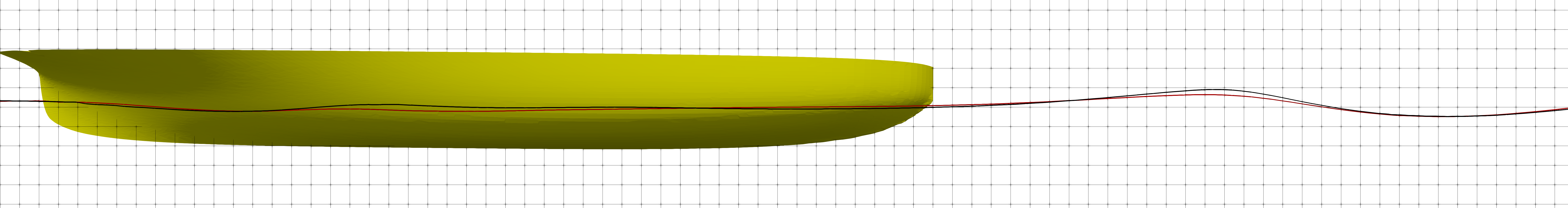}
\caption{Comparison between wave pattern generated by the PHF and optimal hulls, computed by using the RANSE solver
        {\tt interFoam} from the suite {\tt OpenFOAM}\textcopyright \cite{OpenFOAM}. On top, is the top view of the
        wave patterns, where the optimal hull is presented in the lower part. On the bottom, is the comparison of a
        longitudinal wave cut of the PHF (black) and optimal (red) hulls. The trace of the cut is reported as a black
        line in the top view.
         }\label{viste}
\end{center}
\end{figure}

In figure \ref{etaBEM} and \ref{viste} we have a snapshot of the wave pattern produced by the PHF and the optimal hull.
Figure \ref{etaBEM} reports the differences between the wave profile on the hull (and on the centerplane) as simulated
by the BEM solver. In the picture, the hull is located in between $x=-0.5$ and $x=0.5$, with the bow at $x=-0.5$. The
optimal hull shows a more regular wave pattern along the hull, and the hollow observed in the rear part of the PHF has
disappeared. In the wake, the wave profile is reduced after optimization, which is typically beneficial.

The same conditions have been also simulated by using a RANSE mathematical model, with richer physical content. Results
are substantially confirmed, as previously mentioned. Here we can compare the wave pattern predictions: in figure
\ref{viste} we can see a top view of the wave pattern and a comparison of a longitudinal wave cut, for both the PHF
and optimal configuration. The position of the longitudinal cut is reported as a black line in the upper frame of the
same picture. The reduction in the wave elevation along the hull and in the region close to the stern is clearly evident.
We cannot completely exclude that part of the great success of the optimization activities could be possibly connected
with a slightly inaccurate reproduction of the linesplans of the real PHF shape (both in the drawings reported in
\cite{Zotti2002} and/or in the digitalized ones), so that the performances of the PHF are lower than reality.

\section{Conclusions}

The paper is evidencing the connections between AI and optimization, demonstrating how some techniques classically adopted
in AI can be easily and fruitfully applied as base elements of an optimization algorithm. There are some limitations,
mainly connected with the space dimensionality of the problem: in fact, to consider a large number of design parameters
may imply the requirement of a very large DOE, causing at the same time a huge computational cost for the training of the
MM. This situation may become even harder if the objective function is multimodal so that the number of points required
for the synchronization of the prediction of the MMs during the ML phase becomes also larger than in the present cases.
In these conditions, the use of the classical optimization approach could still represent a more viable solution.

More experiences are needed to better establish the limits of the approach. Also, the use of more than a couple of MM,
to further improve the tuning phase, could be investigated.


\section*{Acknowledgements}

This research was funded by Italian Minister of Instruction, University and Research (MIUR) to support this research with funds coming from PRIN Project 2017 (No. 2017KKJP4X entitled ”Innovative numerical methods for evolutionary partial differential equations and applications”).

\section*{Conflict of interest}

The authors declare that they have no conflict of interest.

\bibliography{Peri-PSI-LM}

\end{document}